\documentclass[11pt]{amsart}
\usepackage{amssymb,amsmath,amsfonts,graphicx,subfigure,pstricks,color,a4wide}
  \usepackage{paralist}
  \usepackage{epsfig} 

\usepackage{cite}

\newtheorem{theorem}{Theorem}[section]

\theoremstyle{definition}

\newtheorem{remark}{Remark}

\begin{document}

\title{Tomographic Reconstruction with Spatially Varying Parameter Selection}

\author{Yiqiu Dong} \thanks{Department of Applied Mathematics and Computer Science, Technical University of Denmark, 2800 Kgs. Lyngby, Denmark ({\tt yido@dtu.dk}).} 
\author{Carola-Bibiane Sch\"onlieb} \thanks{Department of Applied Mathematics and Theoretical Physics, University of Cambridge, Cambridge CB3 0WA, UK ({\tt cbs31@cam.ac.uk}).}

\begin{abstract}
In this paper we propose a new approach for tomographic reconstruction with spatially varying regularization parameter. Our work is based on the SA-TV image restoration model proposed in \cite{SATV} where an automated parameter selection rule for spatially varying parameter has been proposed. Their parameter selection rule, however, only applies if measured imaging data are defined in image domain, e.g. for image denoising and image deblurring problems. By introducing an auxiliary variable in their model we show here that this idea can indeed by extended to general inverse imaging problems such as tomographic reconstruction where measurements are not in image domain. We demonstrate the validity of the proposed approach and its effectiveness for computed tomography reconstruction, delivering reconstruction results that are significantly improved compared the state-of-the-art.
\end{abstract}

\date{\today}

\maketitle


\section{Introduction}

Computed tomography (CT) technique has been involved in many clinical and industrial applications. By sending X-rays through the object of interest, we are able to measure the reduction in intensity on the opposite side due to the attenuation of the X-ray when travelling through the object. Then using reconstruction methods we will obtain a 2D or 3D image of the X-ray attenuation coefficients in the object. 

Since X-rays travel along a straight line, mathematically CT scans have been modelled by line integrals, and the corresponding reconstruction problem becomes to find the function from the knowledge of its line integrals
for which an explicit inversion formula was derived in 
\cite{Natterer:2001}. The assumption of obtaining this analytical solution is to have complete and clean data, i.e., the data includes the line integrals in the continuous setting from all directions and without noise. But in real applications the data are limited, and noise is unavoidable. In order to reconstruct high-quality images from a limited number of noisy measurements, several reconstruction methods are proposed in the literature, see e.g. \cite{Herman, QuintoExterior88, OnLocalTomography, FrikelLimitedAngle2013} as well as the monograph \cite{Natterer:2001} and the many references therein.

In recent years, much attention has been given to variational methods for CT reconstruction \cite{Sidky2008, FrikelSparseRegu2013, Dong2016, Li2016, Dong2018}. Based on the maximum a posterior (MAP) estimator and the assumption of white additive Gaussian noise, the variational model for CT reconstruction can be written in the form
\begin{equation}\label{generalmodel}
\min_{u\in L^{2}(\Omega)} \frac{\alpha}{2}\|A u-f\|^{2}+\mathcal{R}(u),
\end{equation}
where $u\in L^{2}(\Omega)$ is the reconstructed image supported in an open subset $\Omega=(-b,b)^{2}\subset\mathbb{R}^{2}$, $f\in L^{2}([0,\pi)\times(-b,b))$ denotes the CT data (often also called as sinogram), and $A: L^{2}(\Omega)\rightarrow L^{2}([0,\pi)\times(-b,b))$ is regarded as the X-ray transform, which is a bounded linear operator in $L^2(\Omega)$. Furthermore, $\mathcal{R}(u)$ is the regularization term, and $\alpha>0$ is the regularization parameter. 
Typical examples for $\mathcal{R}$ are Tikhonov regularization \cite{Tikhonov}, total variation (TV) \cite{TVmodel}, and several extensions \cite{Weickert1998,Scherzer1998,Bredies2010,Lefkimmiatis2015}. The regularization parameter $\alpha$ in \eqref{generalmodel} controls the trade-off between a good fit to the data and the regularization induced on a minimiser of \eqref{generalmodel} by $\mathcal{R}$. The choice of $\alpha$ is critical for receiving a desirable reconstruction. If $\alpha$ is chosen too large the reconstruction will be under-regularised and might still contain noise and other artefacts due to imperfections in the data. On the other hand, an $\alpha$ that is chosen too small will render an over-regularized solutions in which structural information is lost. How to choose $\alpha$, therefore, is an important and moreover challenging question. Even more so, as is the case in \eqref{generalmodel}, when the regularization term is defined in image domain but the data-fitting term is defined in measurement space. A classic approach for regularization parameter selection in this case is Morozov's discrepancy principle \cite{Morozov}, which has been used for CT reconstruction, e.g. in \cite{Hansen:2010, CTpara2015}.  
 
In this paper, we propose a variational method with spatially varying regularization parameter for CT reconstruction. The idea behind the spatially varying regularization parameter is that objects imaged with CT usually contain structures of different scales. In order to preserve these scales in the CT image reconstruction while still removing noise and other artefacts, different regularization parameters should be assigned to different structural scales in the reconstruction. Our proposed method is based on the work in \cite{SATV} for image denoising with spatially adaptive total variation regularisation. Therein, the structural scales are defined in image domain via a spatially varying regularisation parameter that is built into the data fitting term.  Since in CT the data fitting term in \eqref{generalmodel} is defined in measurement space it is not immediately clear how this approach can be applied to this case, or more generally to inverse imaging problems in which reconstruction (image) space and measurement space do not coincide. The variational model that we propose here circumvents this problem by introducing a new variable in \eqref{generalmodel} to split the tomographic reconstruction step and the spatially varying parameter estimation and regularisation step. This new variational model is introduced in \eqref{newmodel} in Section \ref{sec:newmodel}. By using the spatially varying regularization parameter, we will demonstrate that the new method provides much better CT reconstruction results compared to total variation regularisation with scalar $\alpha$. Moreover, the proposed model can be used for reconstruction problems for other inverse imaging problems.      
 
Our paper is organized as follows. In Section \ref{sec:review} we review the work on spatially varying regularization parameter selection in \cite{SATV}. Section \ref{sec:newmodel} describes our new variational method for CT reconstruction in detail. In Section \ref{sec:numerics} we present numerical results from simulated and real measured data. Finally, conclusions are drawn in Section \ref{sec:conclusion}.

\section{Review of the SA-TV Method}\label{sec:review}

Since objects are usually comprised of structures with different scales, locally varying regularization is desirable. In \cite{SATV} a fully automated adjustment strategy for a spatially varying regularization parameter for image denoising and deblurring was proposed that is based on local variance estimators. The resulting regularization method is called spatially adaptive total variation (SA-TV) method, and the corresponding variational model is
\begin{equation}\label{SATV}
\min_{u\in BV(\Omega)} \frac{1}{2}\int_{\Omega} \lambda(x)|K u(x)-z(x)|^{2}\ dx+ \int_{\Omega} |Du|,
\end{equation}
where $K\in\mathcal{L}(L^{2}(\Omega))$ is a blurring operator, $z\in L^{2}(\Omega)$ is a blurred image corrupted by additive white Gaussian noise with mean $0$ and variance $\sigma^{2}$, $\lambda(x)$ is in $L^{\infty}(\Omega)$ and bounded by $[\varepsilon, \bar{\lambda}]$ with $\varepsilon>0$, and $BV(\Omega)$ is the space of functions of bounded variation. Here, $u\in BV(\Omega)$ if $u\in L^1(\Omega)$ and its total variation (TV)
\[
\int_\Omega |Du|= \sup\bigg\{\int_\Omega u\;\mathrm{div} \vec{v} \;
dx:\vec{v}\in(C_0^\infty(\Omega))^2,\|\vec{v}\|_{\infty}\leq1\bigg\}
\]
is finite, where $(C_0^\infty(\Omega))^2$ is the space of vector-valued functions with compact support in $\Omega$. The space $BV(\Omega)$ endowed with the norm $\|u\|_{BV(\Omega)}=\|u\|_{L^1(\Omega)}+\int_\Omega |Du|$
is a Banach space; see, e.g., \cite{BVBook2}. The capability of the SA-TV method strongly depends on the correct selection of the parameter function $\lambda$. 

In order to obtain a locally varying $\lambda$, the idea behind the SA-TV method is to find $\lambda(x_{0})$ for all $x_{0}\in\Omega$ such that the corresponding restored image $u_\lambda$ satisfies the local constraint
\begin{equation}\label{local1}
\int_{\Omega^{\omega}_{x_{0}}} |K u_\lambda(x)-z(x)|^{2}\ dx\leq\sigma^{2}|\Omega^{\omega}_{x_{0}}|,
\end{equation}
where $\Omega^{\omega}_{x_{0}}$ denotes a subset of $\Omega$ with size $[-\frac\omega2,\frac\omega2]\times[-\frac\omega2,\frac\omega2]$ and centered at $x_{0}$, and $|\Omega^{\omega}_{x_{0}}|$ gives its area. Roughly, the constraint \eqref{local1} means that in each local region $\Omega^{\omega}_{x_{0}}$ we expect that the local variance of the residual is less than the noise variance, which can be understood as claiming that for a correctly chosen $\lambda$ nearly only noise is left in the residual and consequently the selected $\lambda(x_{0})$ would automatically depend on the noise level and the scale of textures in this region. For example, if $z$ is rather homogeneous in $\Omega^{\omega}_{x_{0}}$, then we expect that the constraint \eqref{local1} would be satisfied for a small $\lambda$; on the other hand, if $z$ features a lot of small scale textures inside $\Omega^{\omega}_{x_{0}}$, a larger $\lambda$ is needed to preserve these textures in the restored image and only leave the noise in the residual. 

Because the decision on the acceptance or rejection of a local parameter value relies on the scale of textures in the local region, the method potentially requires that the residual is defined in image domain and the textures in $z$ can be easily distinguished, i.e., the operator $K$ is limited to small transformations in image space, e.g. a blurring operator that performs only slight blurring. Due to this limitation, the SA-TV method cannot be directly applied to general inverse problems where $z$ is not in image domain.

\section{Tomographic Reconstruction Method with Spatially Varying $\lambda$}\label{sec:newmodel}

In this section, we follow the same adjustment strategy for a spatially varying regularization parameter as proposed in \cite{SATV} and presented in Section \ref{sec:review} but propose a novel extension of their approach that can be applied to general inverse imaging problems. We exemplify our proposed scheme for the problem of CT reconstruction. 

Considering the TV regularization in the CT reconstruction model \eqref{generalmodel}, we obtain
\begin{equation}\label{CTmodel}
\min_{u\in BV^{+}(\Omega)} \frac{\alpha}{2}\|A u-f\|^{2}+ \int_{\Omega} |Du|,
\end{equation}
where $BV^{+}(\Omega)=\{u\in BV(\Omega): u(x)\geq 0\}$. In \eqref{CTmodel} the regularization is based on the smoothness assumption on the reconstructed image $u$, and the data-fitting term arises from the CT forward model 
\[
f=Au+\epsilon
\]
with the assumption that the noise $\epsilon$ follows a Gaussian distribution with mean 0 and variance $\sigma^2$. Since $A$ denotes X-ray transform, i.e. a line integral operator, each point in the sinogram $f$ gives a value of a line integral, which is global information, such that $f$ is very smooth and carries different singularities comparing to $u$, see the detailed study in \cite{QuintoMLA} by applying microlocal analysis. Hence, we cannot select the regularization parameter in the same way as in \eqref{local1} according to the scale of textures in the residual. 

To circumvent this problem, we introduce a new variable $w$ and split the data-fitting term and regularization term in \eqref{CTmodel}. The new variational model that we propose is as follows:
\begin{equation}\label{newmodel}
\min_{u, w\in BV^{+}(\Omega)} \mathcal{J}(u,w):=\frac{\alpha}{2}\|A w-f\|^{2}+\frac12\int_{\Omega} \lambda(x)(w(x)-u(x))^{2}\ dx+ \int_{\Omega} |Du|,
\end{equation}
where $\alpha>0$ and $0<\varepsilon\leq\lambda(x)\leq\bar{\lambda}$. Because the sinogram $f$ shows global information on the objects, it is difficult to define a local constraint for choosing a locally different $\alpha$. With the extra variable $w$, however, the new added quadratic term is defined in image domain and $w$ can take over the former role of $z$ in \eqref{local1}. We therefore can introduce a spatially varying parameter based on local image textures in $w$. Note that in the new model the residual $w-u$ generally carries the artifacts from the CT reconstruction due to the noise in the sinogram and the image textures. 

Since $\lambda(x)\geq\varepsilon>0$, the model is strictly convex. Hence, existence and uniqueness of solutions to \eqref{newmodel} is a straightforward exercise, whose result is summarised in the following theorem.
\begin{theorem}
Let $f$ be in $L^{2}([0,\pi)\times(-b,b))$, $A: L^{2}(\Omega)\rightarrow\in L^{2}([0,\pi)\times(-b,b))$ be a bounded linear operator, and $\lambda\in L^{\infty}(\Omega)$ be bounded in $[\varepsilon, \bar{\lambda}]$. Then, the model \eqref{newmodel} exists a unique solution.
\end{theorem}
    
To solve the minimization problem in \eqref{newmodel}, we use an alternating optimisation algorithm, which starts from an initial guess $w^{0}\in BV^{+}(\Omega)$ and follows an iterative scheme
\begin{align}
u^{k+1}=\underset{u\in BV^{+}(\Omega)}{\operatorname{argmin}}\  & \frac12\int_{\Omega} \lambda(x)(w^k(x)-u(x))^{2}\ dx+ \int_{\Omega} |Du|,\label{u-problem}\\
w^{k+1}=\underset{w\in BV^{+}(\Omega)}{\operatorname{argmin}}\ &  \frac{\alpha}{2}\|A w-f\|^{2}+\frac12\int_{\Omega} \lambda(x)(w(x)-u^{k+1}(x))^{2}\ dx,\label{w-problem}
\end{align}
which is solved till numerical convergence, i.e. the difference of iterates is smaller than a prescribed tolerance.

In \eqref{u-problem}, the solution satisfies the minimum-maximum principle, i.e., $\inf_{\Omega}w^k\leq u^{k+1}\leq\sup_{\Omega} w^k$, and we know that $w^k(x)\geq0$ according to the constraint in \eqref{w-problem}, then the $u$-subproblem is equivalent to
\begin{equation}
u^{k+1}=\underset{u\in BV(\Omega)}{\operatorname{argmin}} \frac12\int_{\Omega} \lambda(x)(w^k(x)-u(x))^{2}\ dx+ \int_{\Omega} |Du|,\label{u-problem2}
\end{equation}
which is identical to the model \eqref{SATV} in the SA-TV method with $K$ being the identity operator. Hence, we can apply the SA-TV method proposed in \cite{SATV} to solve the $u$-subproblem with automatically adjusted regularization parameter $\lambda$. The $w$-subproblem in \eqref{w-problem} is a least-squares problem, and we can solve it efficiently by using e.g. the CGLS method \cite{CGLS} followed by a non-negativity projection. The selection of the parameter $\alpha$ is done by using the discrepancy principle. In our method, $\lambda(x)$ is only updated in the first $k_0$ iterations, after that we keep $\lambda$ fixed and let $(u^{k+1}, w^{k+1})$ converge. In all numerical tests, we set $k_0=5$. 

With a fixed $\lambda$, the convergence of the sequence $\{(u^{k}, w^k)\}$ is guaranteed according to 
\[
\mathcal{J}(u^{k+1},w^{k+1})\leq \mathcal{J}(u^{k+1},w^{k})\leq \mathcal{J}(u^{k},w^{k})
\]
as well as the fact that $\mathcal{J}(u, w)$ is bounded below by zero. Then, taking the strict convexity of the objective function $\mathcal{J}$ into account, we obtain the convergence of the algorithm.

\begin{theorem}
Assume that $w^{0}\in BV^{+}(\Omega)$ and $\lambda\in L^\infty(\Omega)$ is fixed and bounded in $[\varepsilon, \bar{\lambda}]$. Then the sequence $\{(u^k, w^k)\}$ converges to the unique minimizer of the problem \eqref{newmodel}. 
\end{theorem}

\begin{remark}
Note that due to the $\lambda$-term in \eqref{w-problem}, the proposed method is different to a two-stage method, which reconstructs image $w$ first by solving a least-squares problem $\min_{w}\|Aw-f\|^{2}$ and then uses the SA-TV method to post-process the image. In our approach the variables $u$ and $w$ are correlated through the $\lambda$-term in the model.
\end{remark}
\begin{remark}
In our method, the spatially varying parameter is determined according to the ``noise'' variance in the residual and the scale of textures in local regions, which is the same as in the SA-TV method. Hence, similar as in \cite{SATGV} the new method also can be extended to apply other regularization terms.
\end{remark}

\section{Numerical Results}\label{sec:numerics}

In this section we provide numerical results for simulated as well as real data to study the behaviour of the proposed reconstruction method with spatially varying regularization parameter. In our method, the parameter $\alpha$ is chosen by applying the discrepancy principle with given or estimated noise variance $\sigma^2$. We apply the SA-TV method to solve the $u$-subproblem \eqref{u-problem}, using the default setting suggested in \cite{SATV} where the only input is the ``noise'' level in the reconstruction, which for the real data is estimated by calculating the variance of the dark background. The stopping criteria in our method is 
\[
\frac{\|u^{k+1}-u^k\|}{\|u^k\|}\leq 10^{-4}.
\]

\begin{figure}[t]
\begin{center}
\begin{minipage}[t]{4.7cm}
\includegraphics[height=4.7cm]{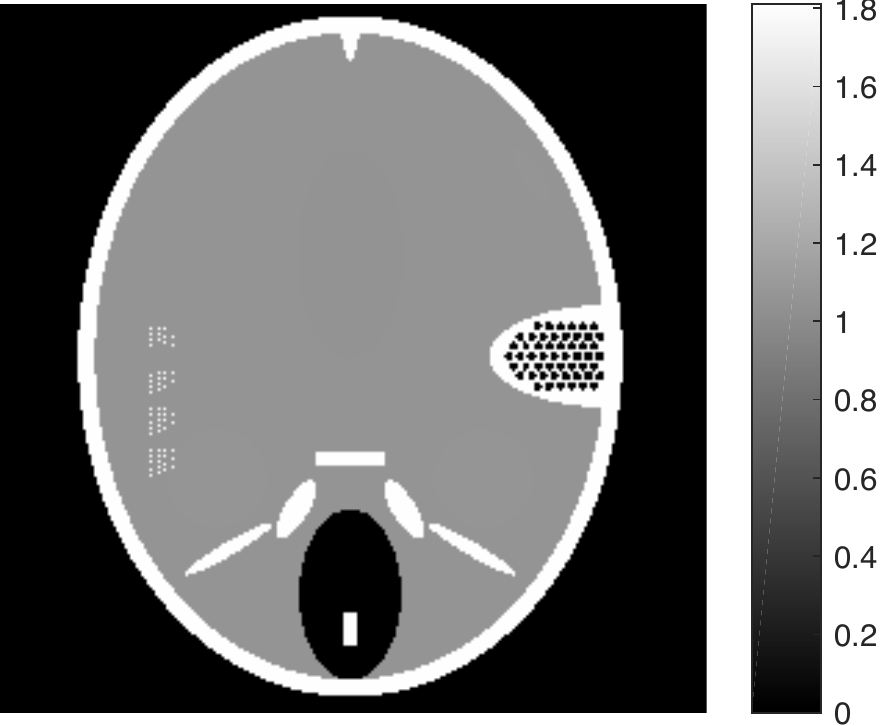}\\
\centering{(a)}
\end{minipage}\hspace{19mm}
\begin{minipage}[t]{4.7cm}
\includegraphics[height=4.6cm]{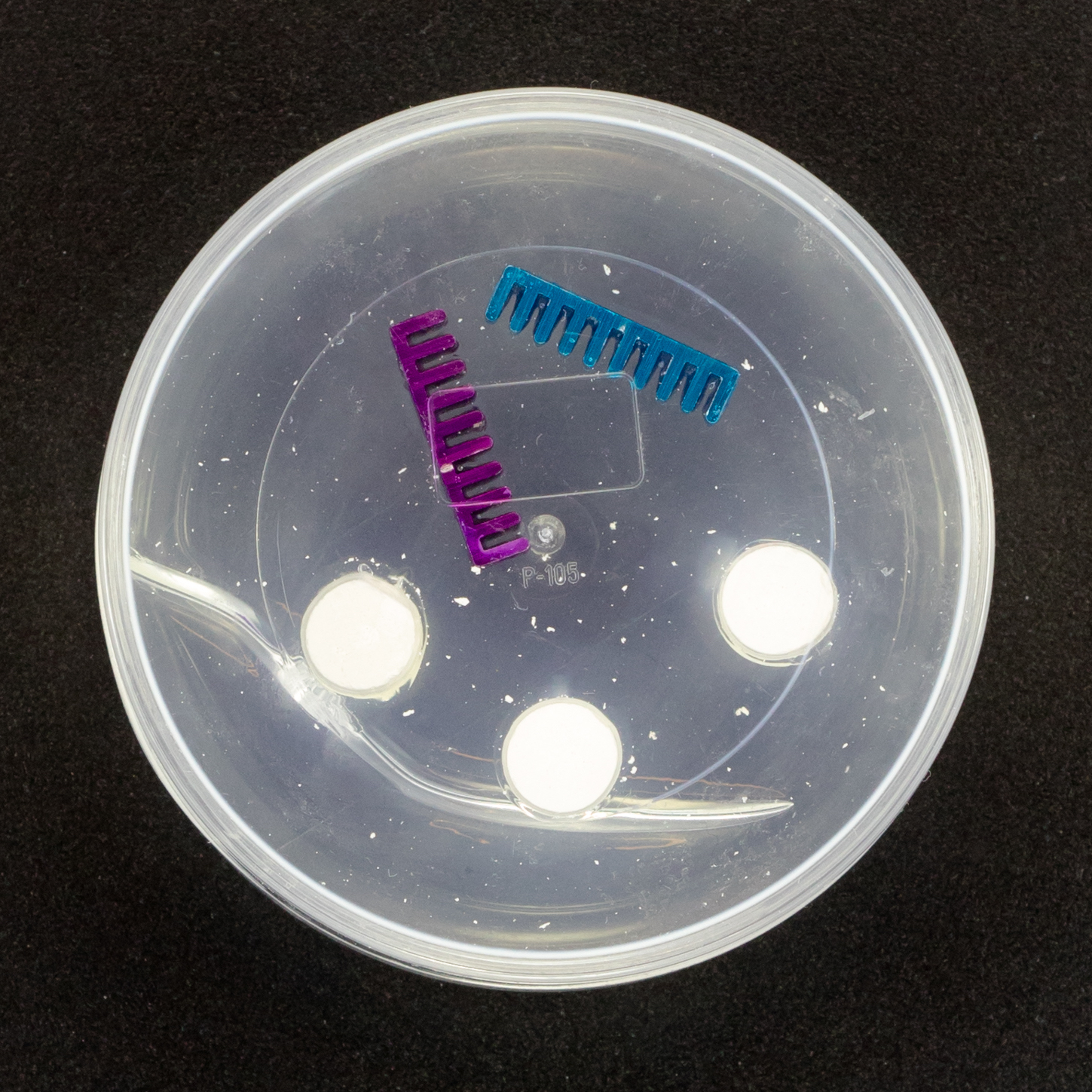}\\
\centering{(b)}
\end{minipage}
\end{center}
\caption{\label{phantom} Phantoms for tests. (a) Head phantom used for simulation \cite{headphantom}, (b) gel phantom used for real X-ray scan.}
\end{figure}

For the experiments, we test our method on a simulated head phantom from \cite{headphantom} and a real gel phantom shown in Figure \ref{phantom}. Here, we consider a 2D CT scenario. The detector has full coverage of the object at any projection angle, and a constant angular spacing of the rays is set in the interval of $[0,\pi]$.

\begin{figure}[t]
\begin{center}
\begin{minipage}[t]{4.7cm}
\includegraphics[height=4.7cm]{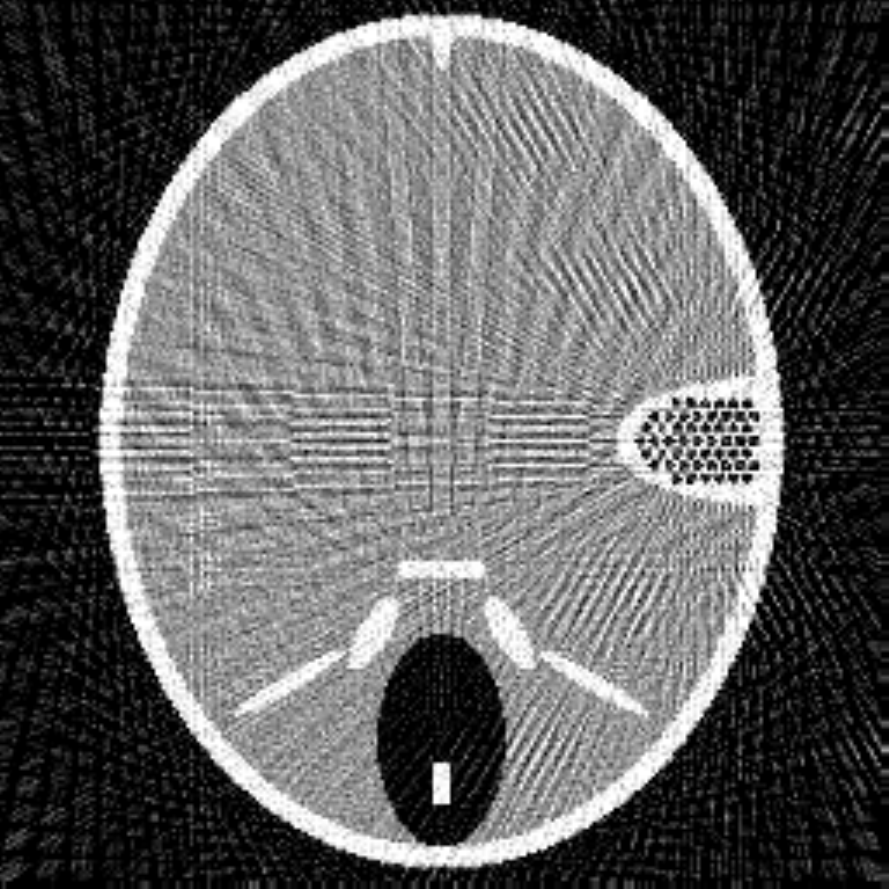}\\
\centering{(a)}
\end{minipage}
\begin{minipage}[t]{4.7cm}
\includegraphics[height=4.7cm]{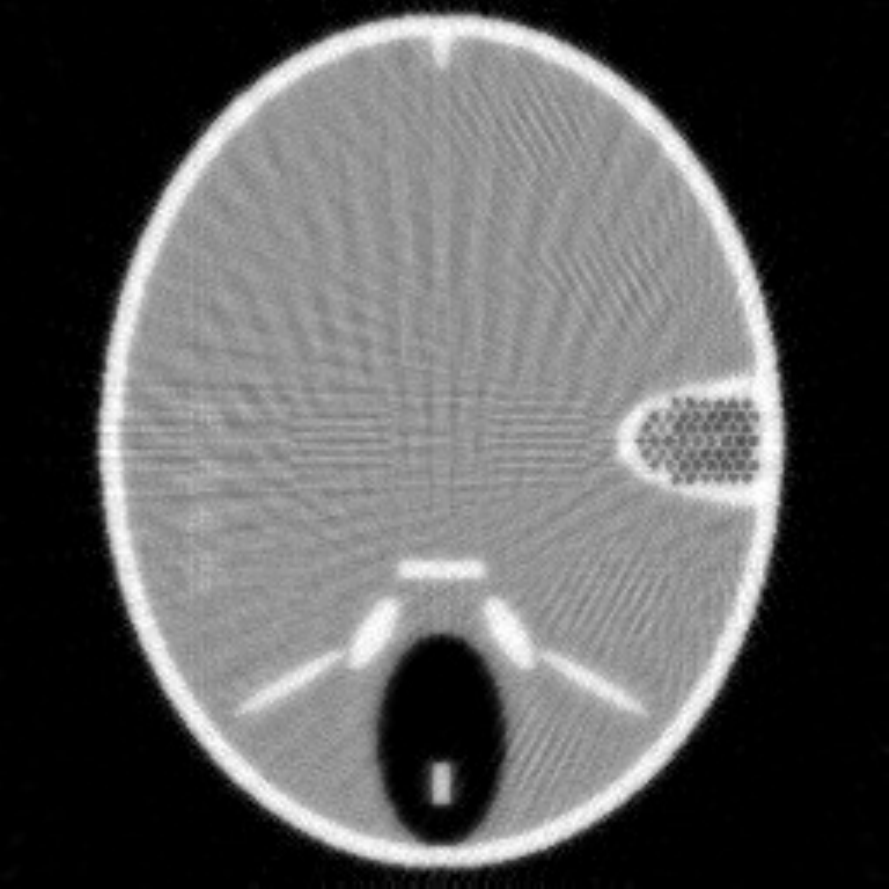}\\
\centering{(b)}
\end{minipage}
\begin{minipage}[t]{4.7cm}
\includegraphics[height=4.7cm]{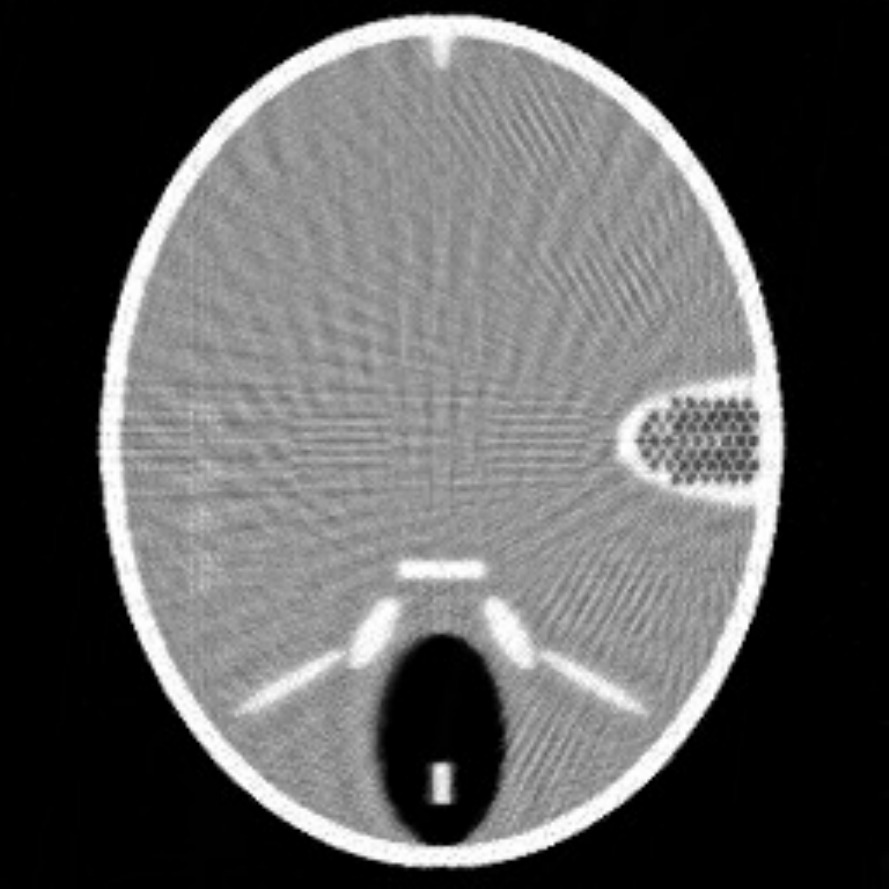}\\
\centering{(c)}
\end{minipage}\\
\begin{minipage}[t]{4.7cm}
\includegraphics[height=4.7cm]{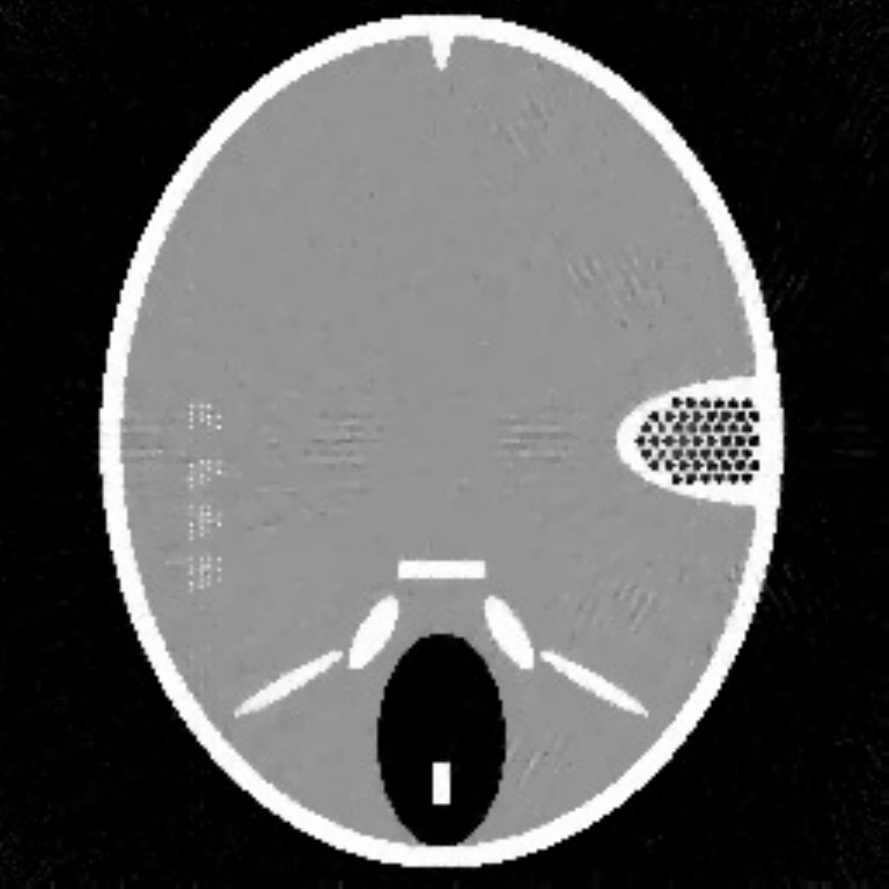}\\
\centering{(d)}
\end{minipage}
\begin{minipage}[t]{4.7cm}
\includegraphics[height=4.7cm]{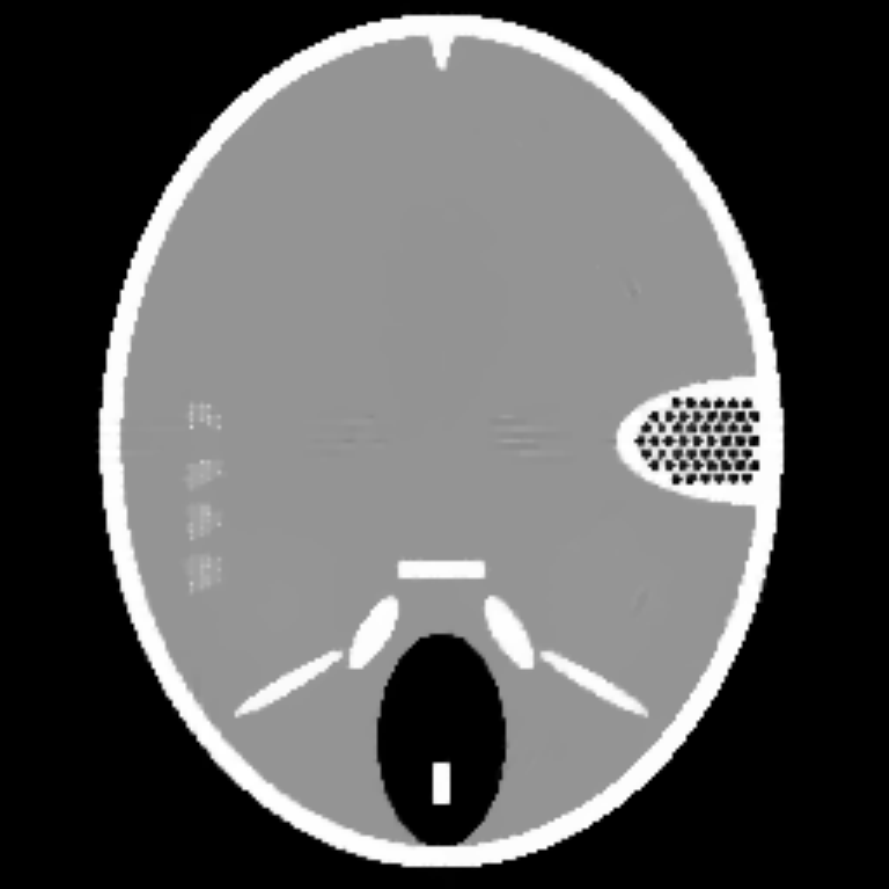}\\
\centering{(e)}
\end{minipage}
\begin{minipage}[t]{4.7cm}
\includegraphics[height=4.7cm]{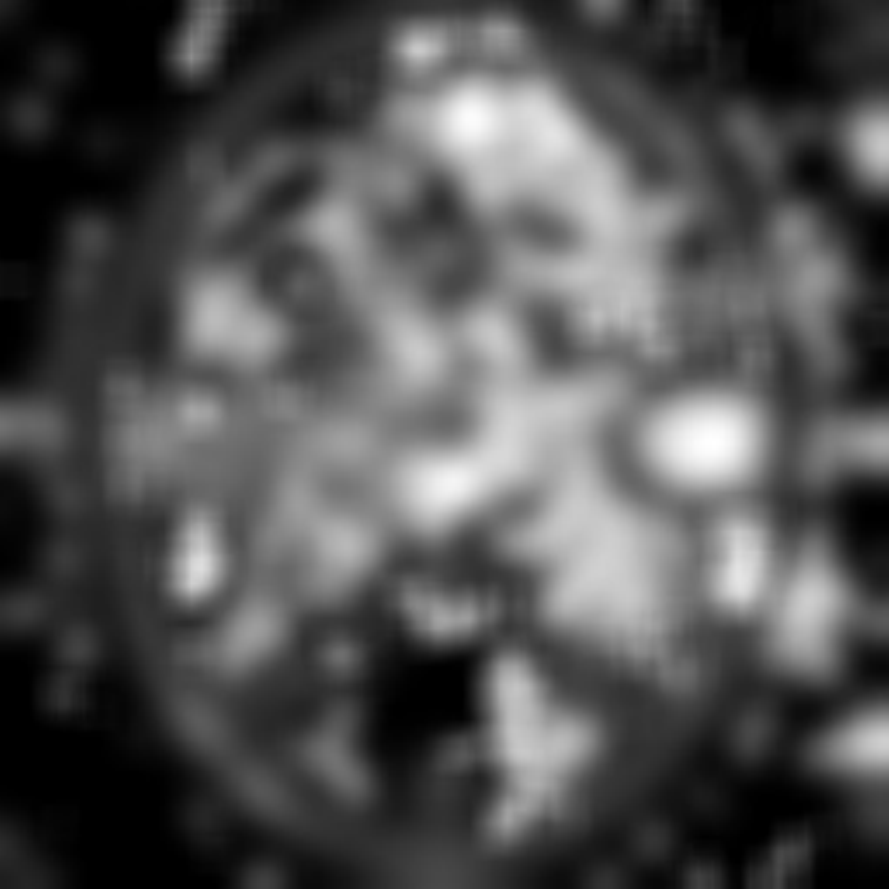}\\
\centering{(f)}
\end{minipage}
\end{center}
\caption{\label{result_head_noise2} Results of different methods for reconstructing the head phantom with underdetermined rate 25\% and relative noise level 0.2. (a) FBP (SNR=0.3253), (b) Landweber (SNR=0.1718), (c) Kaczmarz (SNR=0.1406), (d) L2-TV with scalar $\lambda$ (SNR=0.0644), (e) Our method (SNR=0.0525), (f) $\lambda$ in our method.}
\end{figure}

 \begin{figure}[t]
\begin{center}
\begin{minipage}[t]{4.7cm}
\includegraphics[height=4.7cm]{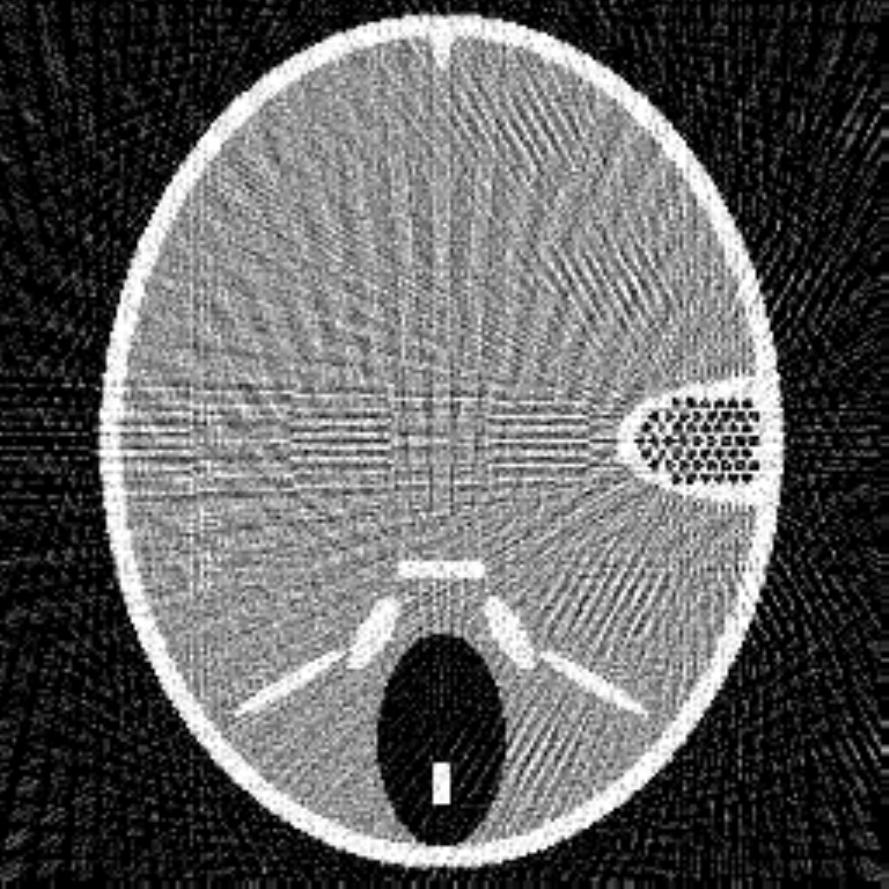}\\
\centering{(a)}
\end{minipage}
\begin{minipage}[t]{4.7cm}
\includegraphics[height=4.7cm]{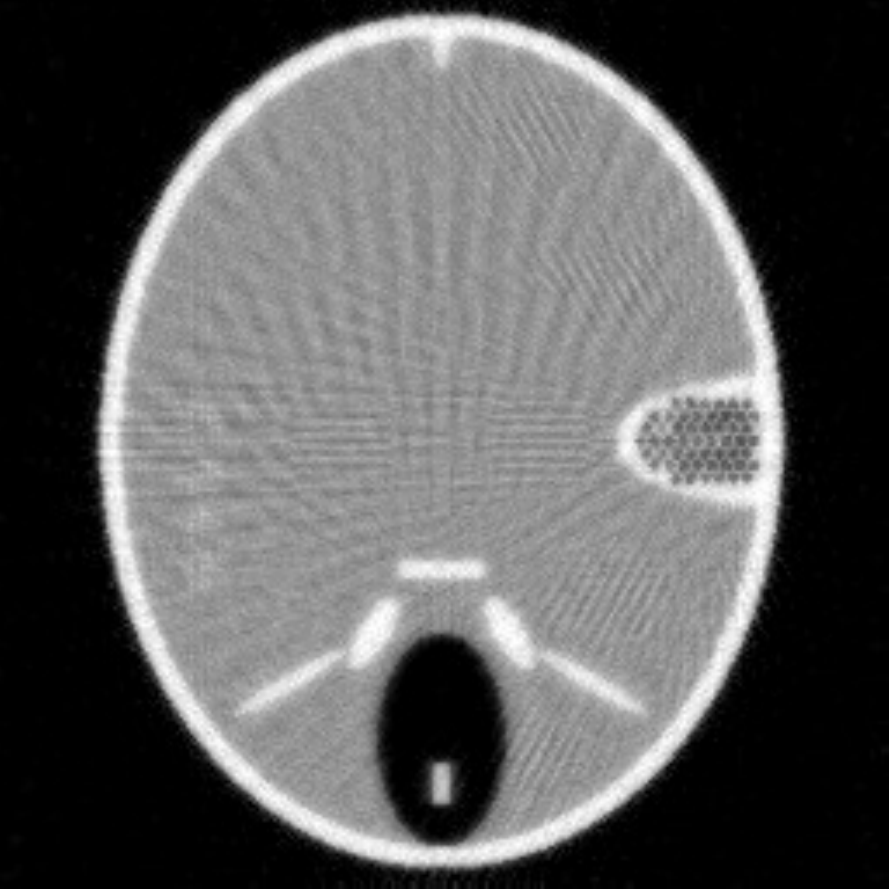}\\
\centering{(b)}
\end{minipage}
\begin{minipage}[t]{4.7cm}
\includegraphics[height=4.7cm]{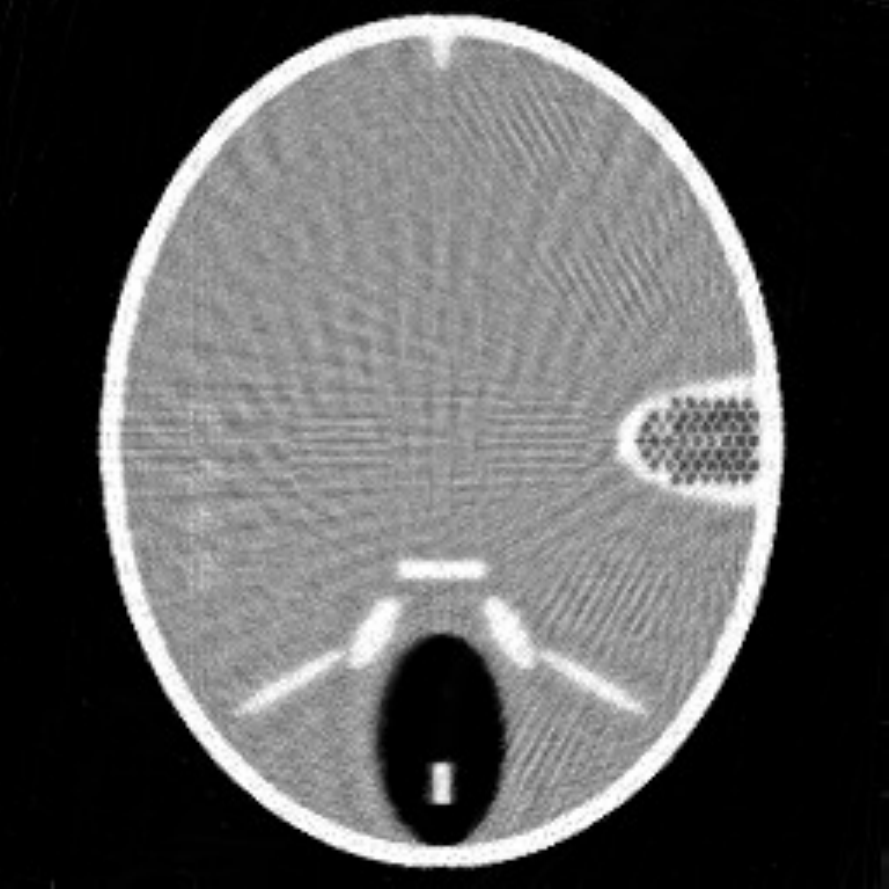}\\
\centering{(c)}
\end{minipage}\\
\begin{minipage}[t]{4.7cm}
\includegraphics[height=4.7cm]{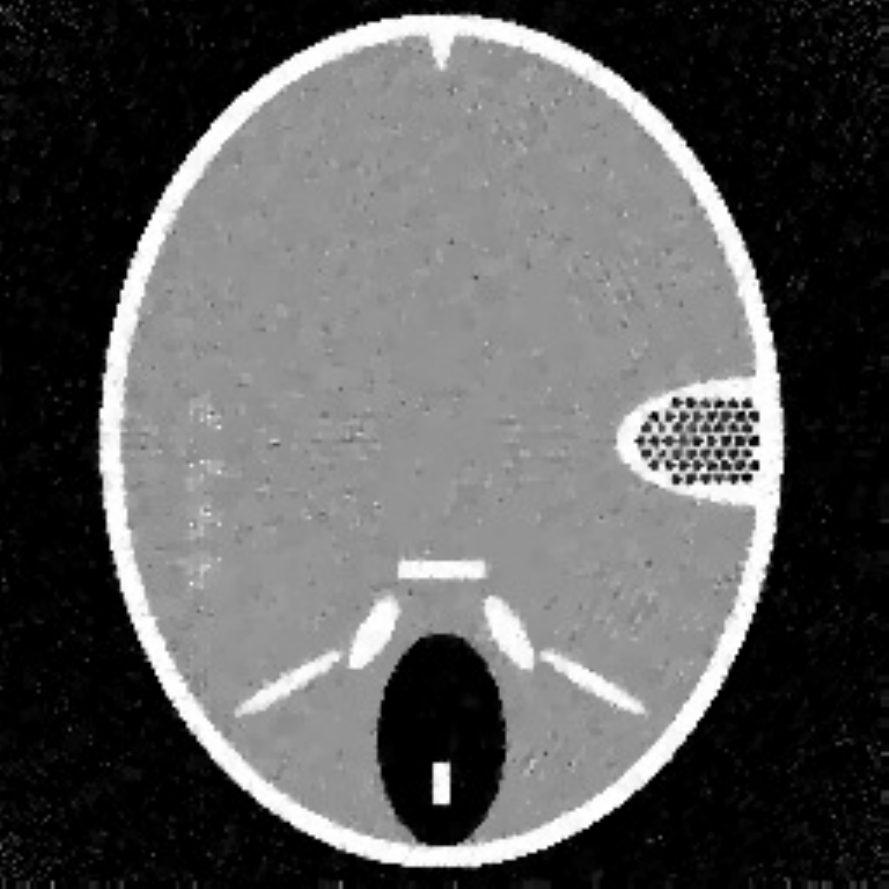}\\
\centering{(d)}
\end{minipage}
\begin{minipage}[t]{4.7cm}
\includegraphics[height=4.7cm]{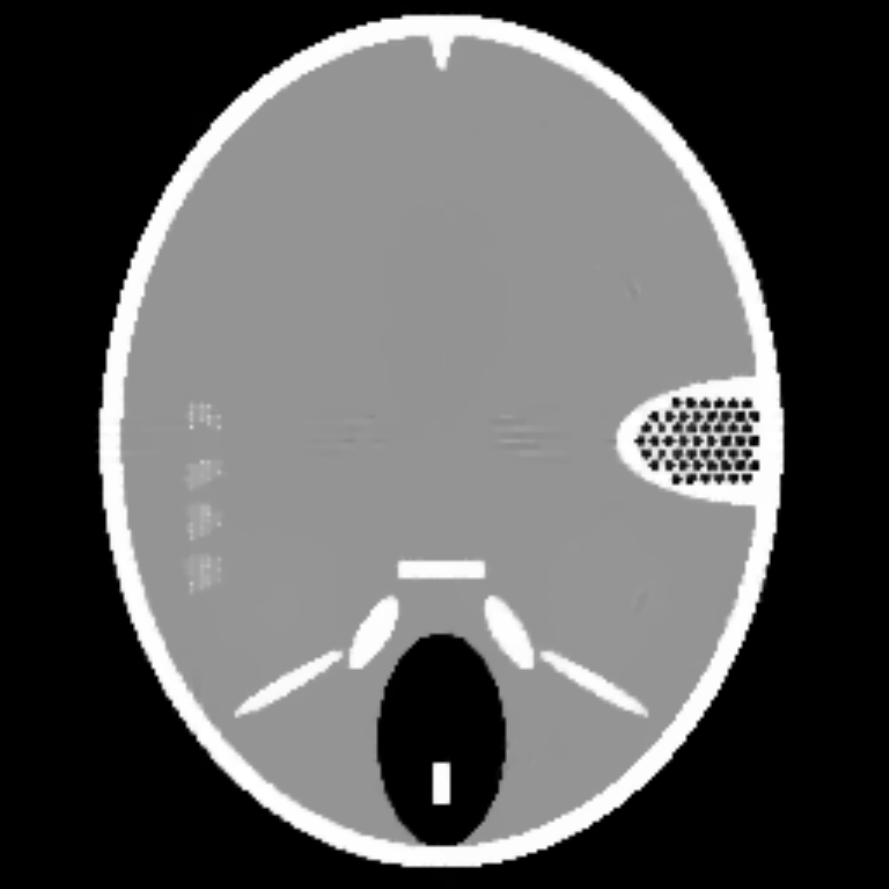}\\
\centering{(e)}
\end{minipage}
\begin{minipage}[t]{4.7cm}
\includegraphics[height=4.7cm]{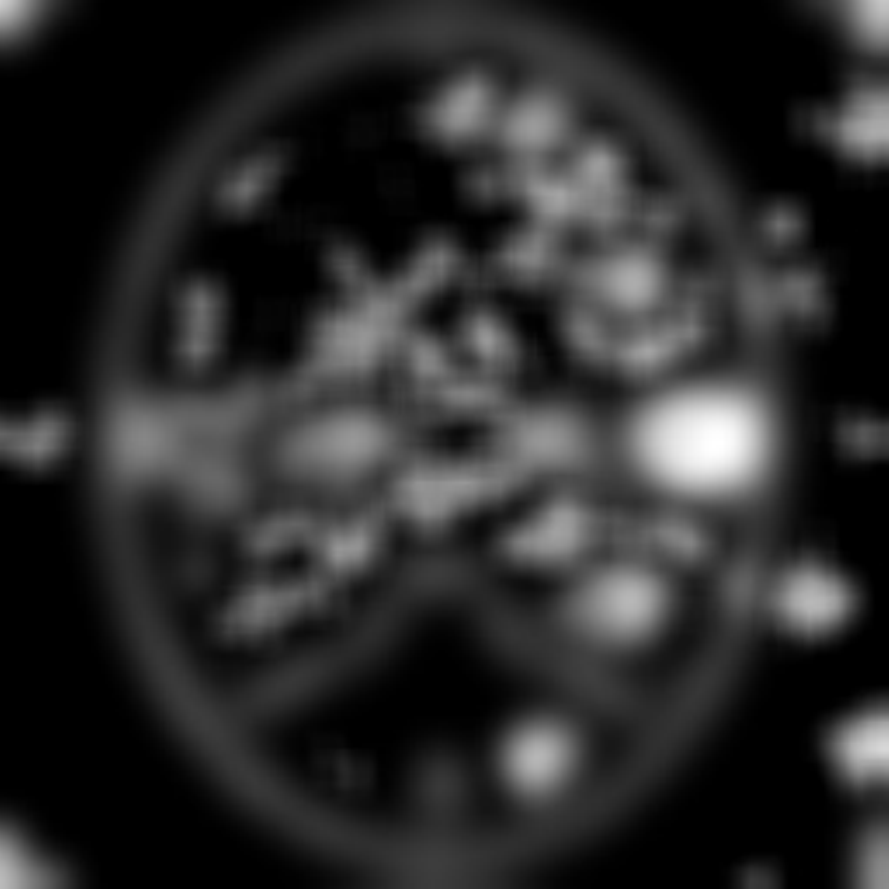}\\
\centering{(f)}
\end{minipage}
\end{center}
\caption{\label{result_head_noise8} Results of different methods for reconstructing the head phantom with underdetermined rate 25\% and relative noise level 0.8. (a) FBP (SNR=0.3505), (b) Landweber (SNR=0.1740), (c) Kaczmarz (SNR=0.1448), (d) L2-TV with scalar $\lambda$ (SNR=0.1072), (e) Our method (SNR=0.0762), (f) $\lambda$ in our method.}
\end{figure}

\textbf{Example 1.} Our first test example is on the simulated head phantom, which is generated in a square domain of $256\times256$ pixels, i.e., there are $256^{2}=65,536$ unknowns. With 362 beams and 45 projection angles, the correspond CT reconstruction problem has an under-determined rate of 25\%. The measurements are given by $f=A\bar{u}+\epsilon$, where $\bar{u}$ is the ground truth (the true attenuation coefficients in the object) and $\epsilon$ denotes the additive white Gaussian noise with the noise level $\|\epsilon\|/\|A\bar{u}\|$.

In order to study the behaviour of our method, we compare it with the filtered back-projection (FBP) algorithm \cite{Natterer:2001}, the Landweber method \cite{Natterer:2001}, the Kaczmarz's method \cite{Herman}, and the L2-TV reconstruction method, which solves the variational model \eqref{generalmodel} with TV regularization as proposed in \cite{TVmodel}. All methods are solved under a non-negativity constraint. Note that in the L2-TV method the regularization parameter $\alpha$ is scalar, which is chosen to give the largest signal-to-noise ratio (SNR).

In Figure \ref{result_head_noise2} and \ref{result_head_noise8}, we give the reconstruction results, which are shown in the same intensity range as the original phantom, from the simulated measurements with the noise level 0.2 and 0.8, respectively. Since the FBP algorithm is according to the analytical formulation of the inverse X-ray transform, it implicitly requires to have continuously measured clean data from the whole 0 to $\pi$ angular range. Therefore, it is not suited for reconstructing from noisy limited data. We can clearly see many stripe artifacts due to the noise and sparse projection angles in the FBP results. Both the Landweber and Kaczmarz's methods perform better than FBP, but there are still some visible artifacts in the reconstruction. By using the TV regularization in the L2-TV and our methods, we potentially assume that the reconstructions are piecewise constant, which evidently reduces the influence of the noise and avoids stripe artifacts. In addition, comparing the results from the L2-TV and our methods, we find that our method suppresses artifacts much better while reconstructing most details. For instance, the grey region in the head and the black dotted region on the right side. With respect to SNR, it is also clear that our method gives the best reconstruction results. In Figure \ref{result_head_noise2} (f) and Figure \ref{result_head_noise8}(f), we also plot the $\lambda$ function obtained by our method. One can see that in the more textured regions $\lambda$ is large in order to preserve the details, and in the more homogenous regions it is small to reduce artifacts.

\begin{figure}[t]
\begin{center}
\begin{minipage}[t]{3.7cm}
\includegraphics[height=3.7cm]{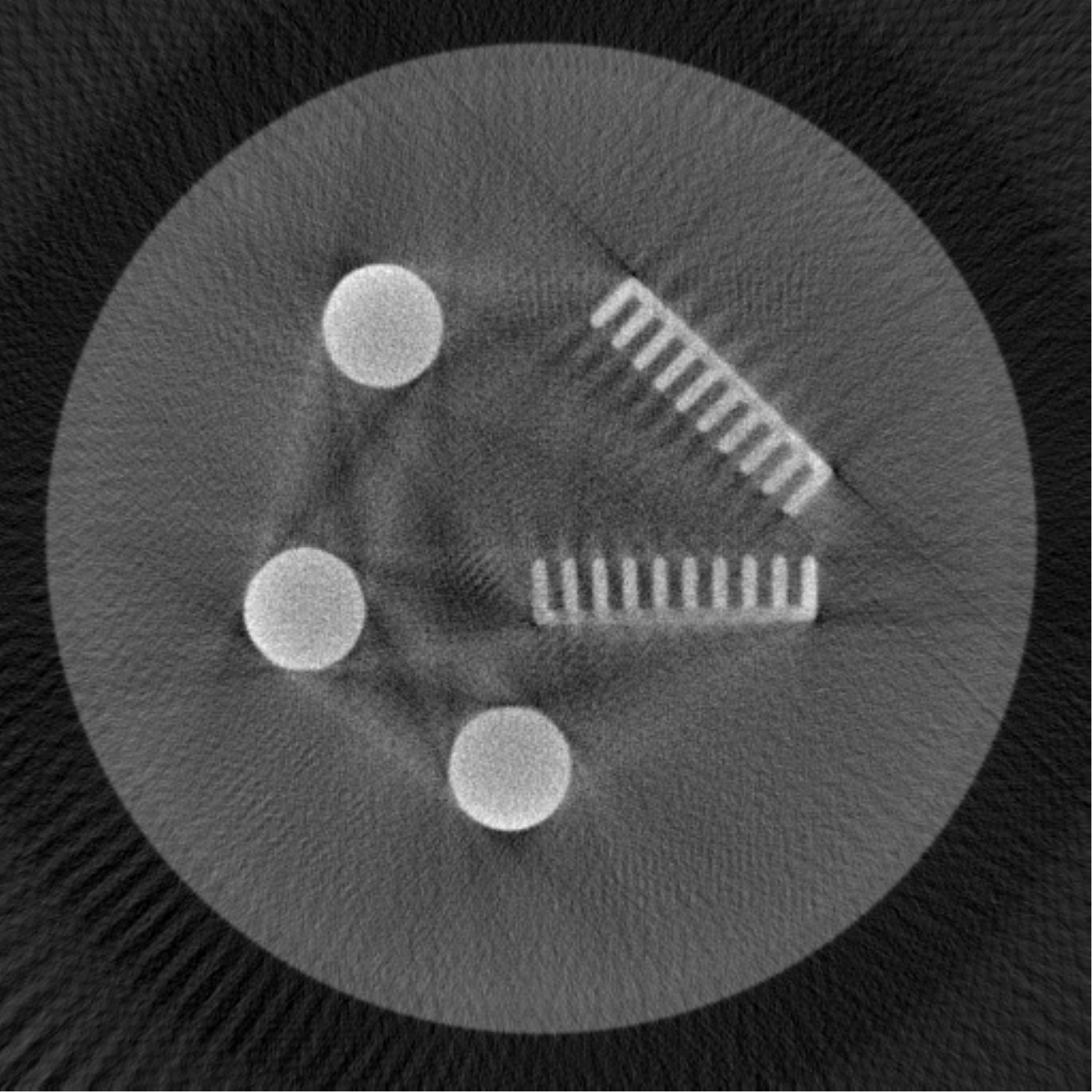}\\
\centering{(a)}
\end{minipage}
\begin{minipage}[t]{3.7cm}
\includegraphics[height=3.7cm]{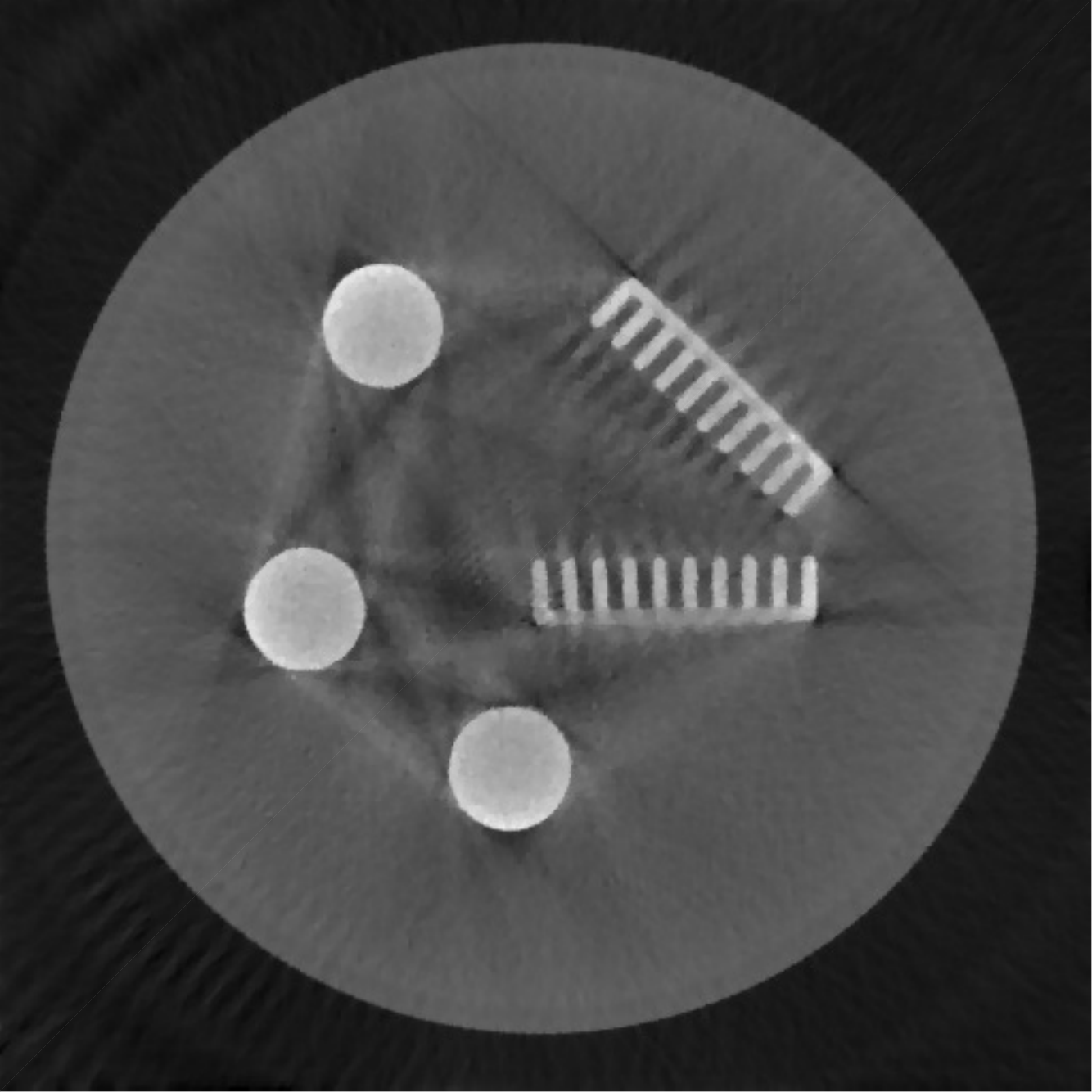}\\
\centering{(b)}
\end{minipage}
\begin{minipage}[t]{3.7cm}
\includegraphics[height=3.7cm]{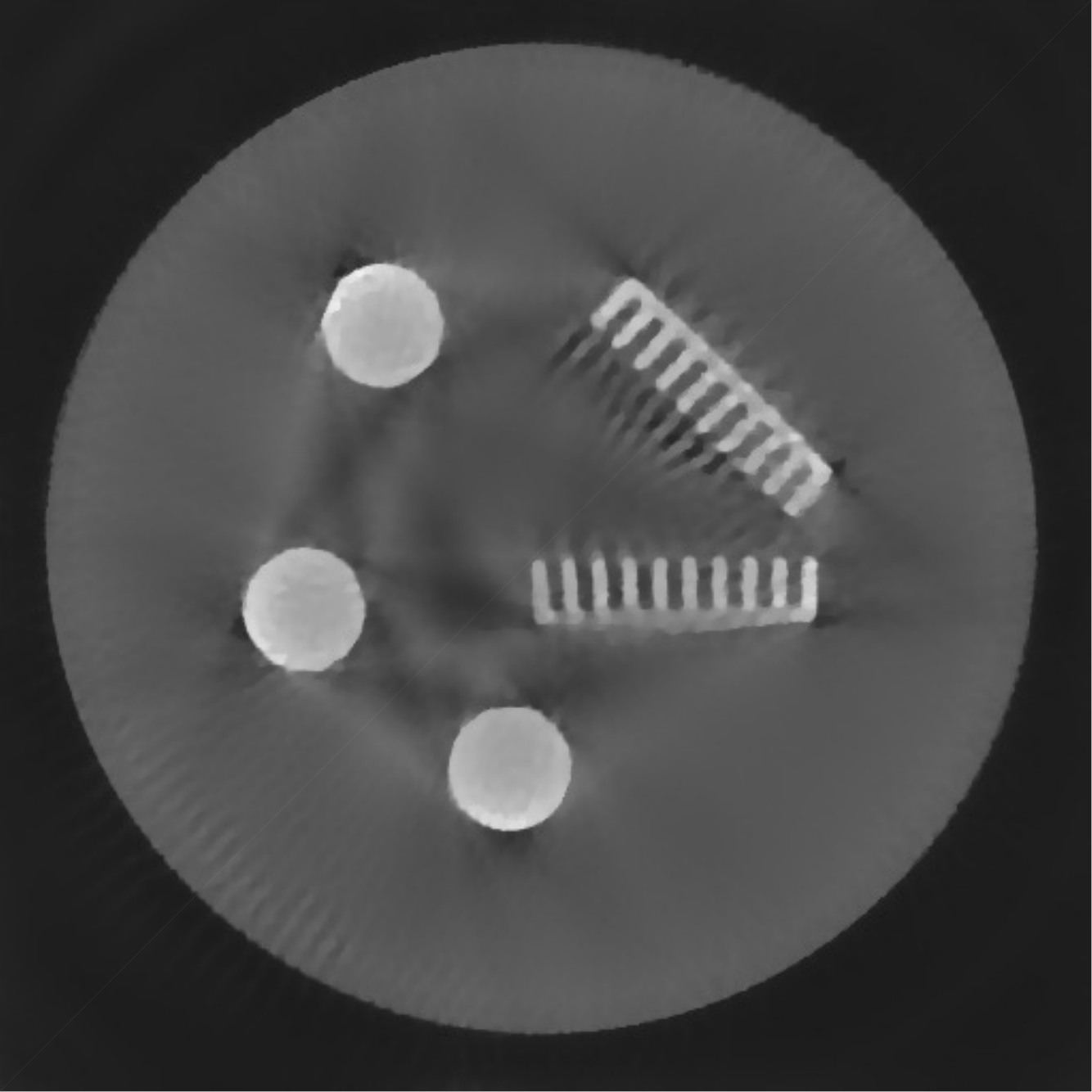}\\
\centering{(c)}
\end{minipage}
\begin{minipage}[t]{3.7cm}
\includegraphics[height=3.7cm]{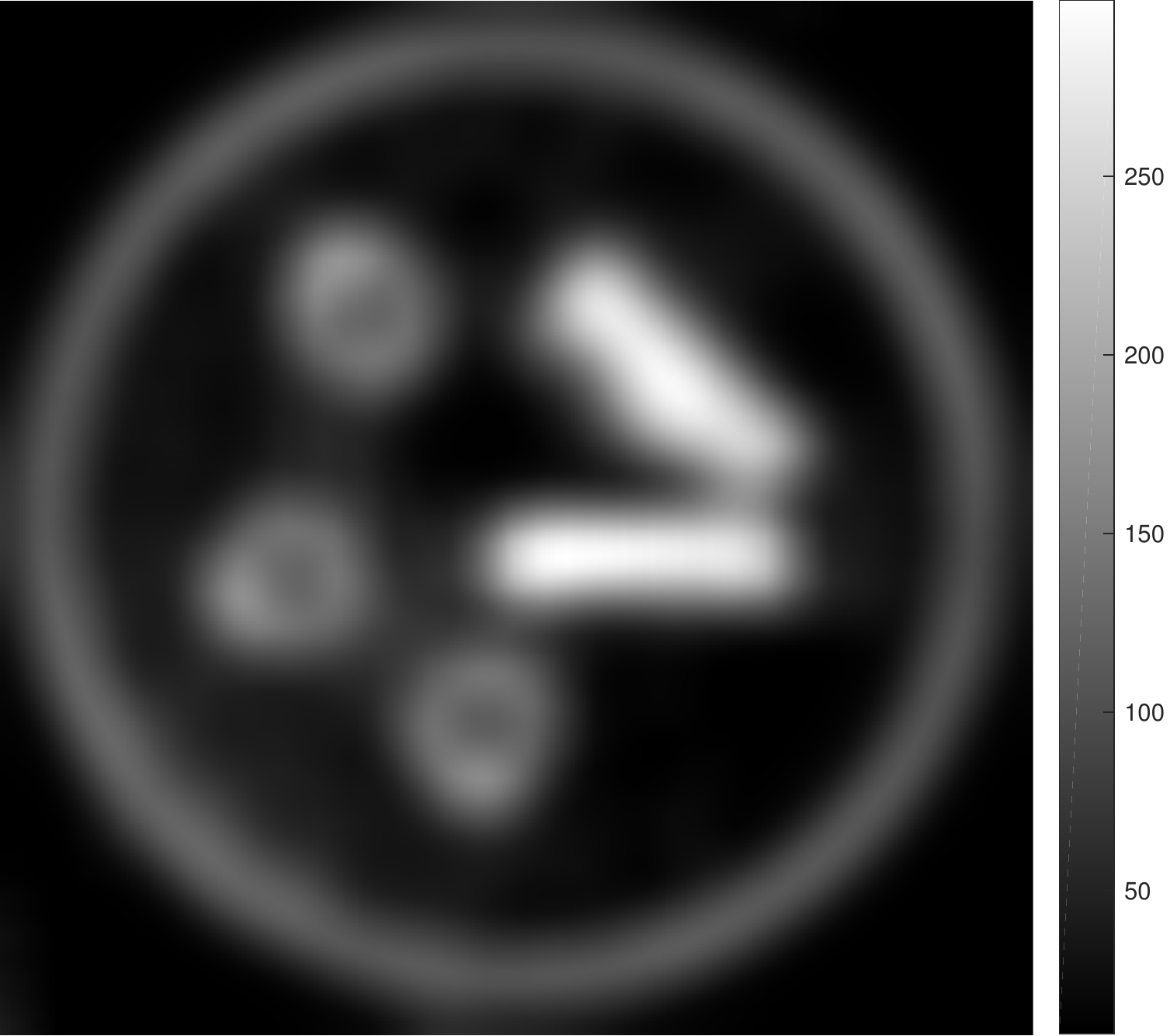}\\
\centering{(d)}
\end{minipage}
\end{center}
\caption{\label{result_gel1} Results of different methods for reconstructing the gel phantom with underdetermined rate 77\%. (a) FBP, (b) L2-TV with scalar $\lambda$, (c) Our method, (d) $\lambda$ in our method.}
\end{figure}

\begin{figure}[t]
\begin{center}
\begin{minipage}[t]{3.7cm}
\includegraphics[height=3.7cm]{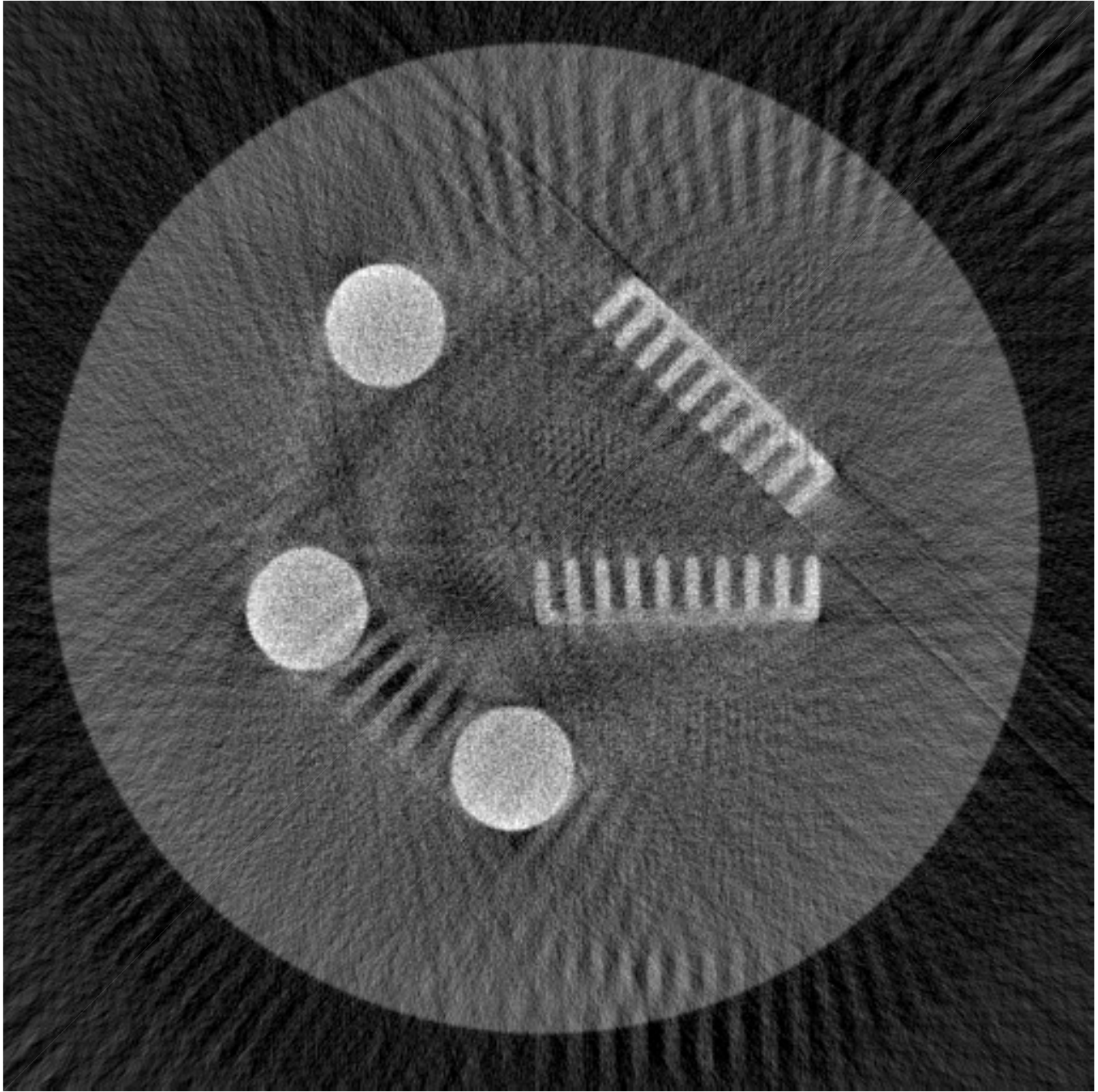}\\
\centering{(a)}
\end{minipage}
\begin{minipage}[t]{3.7cm}
\includegraphics[height=3.7cm]{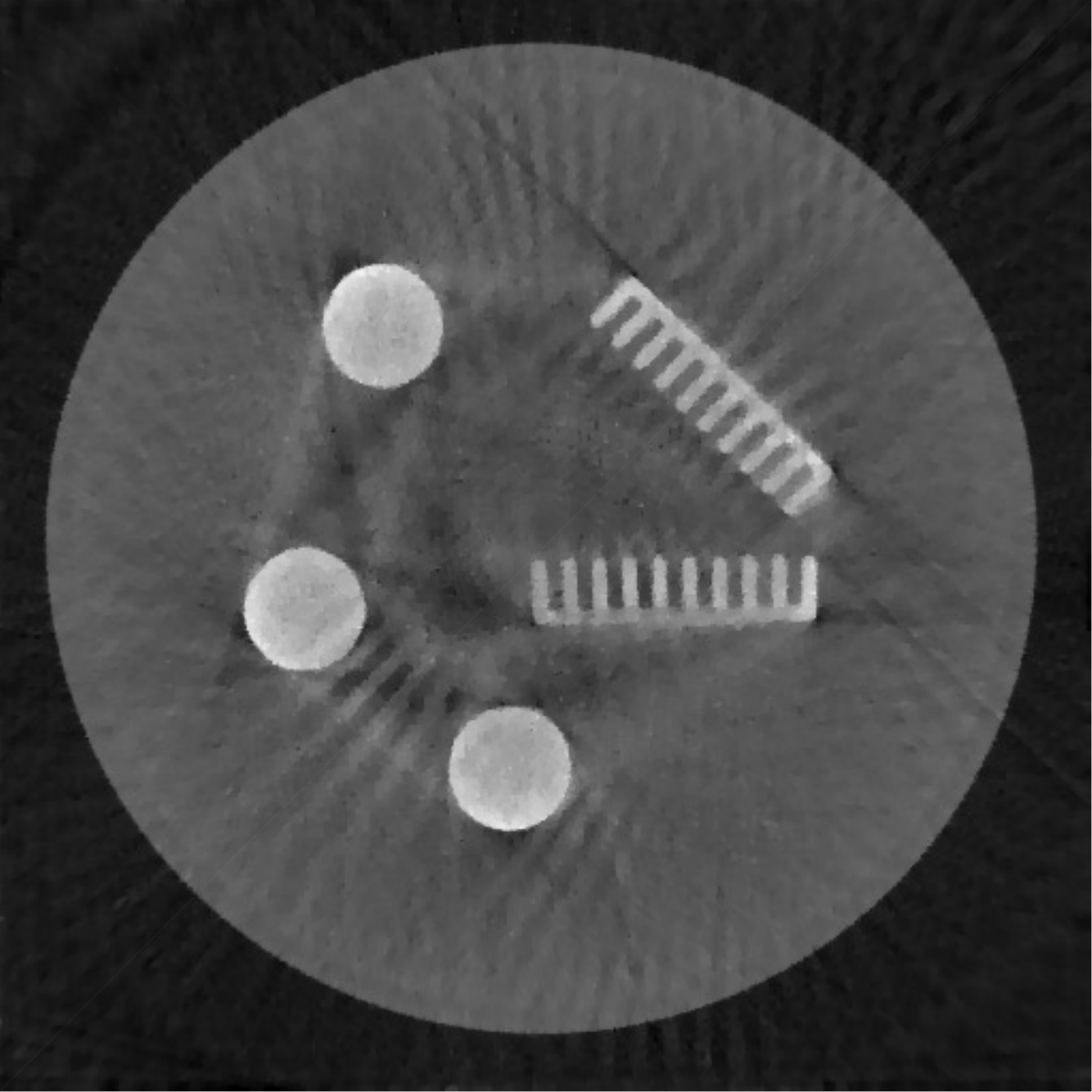}\\
\centering{(b)}
\end{minipage}
\begin{minipage}[t]{3.7cm}
\includegraphics[height=3.7cm]{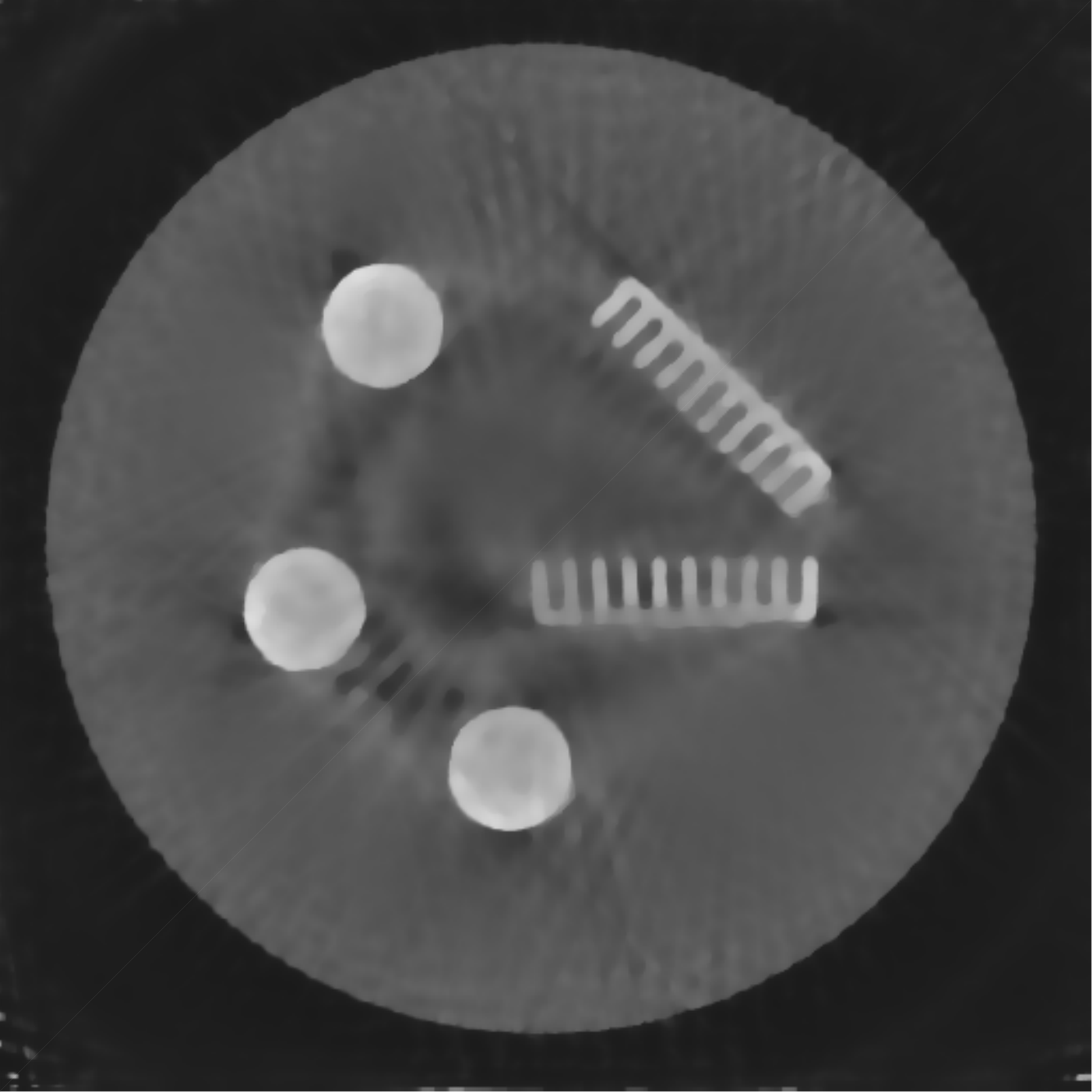}\\
\centering{(c)}
\end{minipage}
\begin{minipage}[t]{3.7cm}
\includegraphics[height=3.7cm]{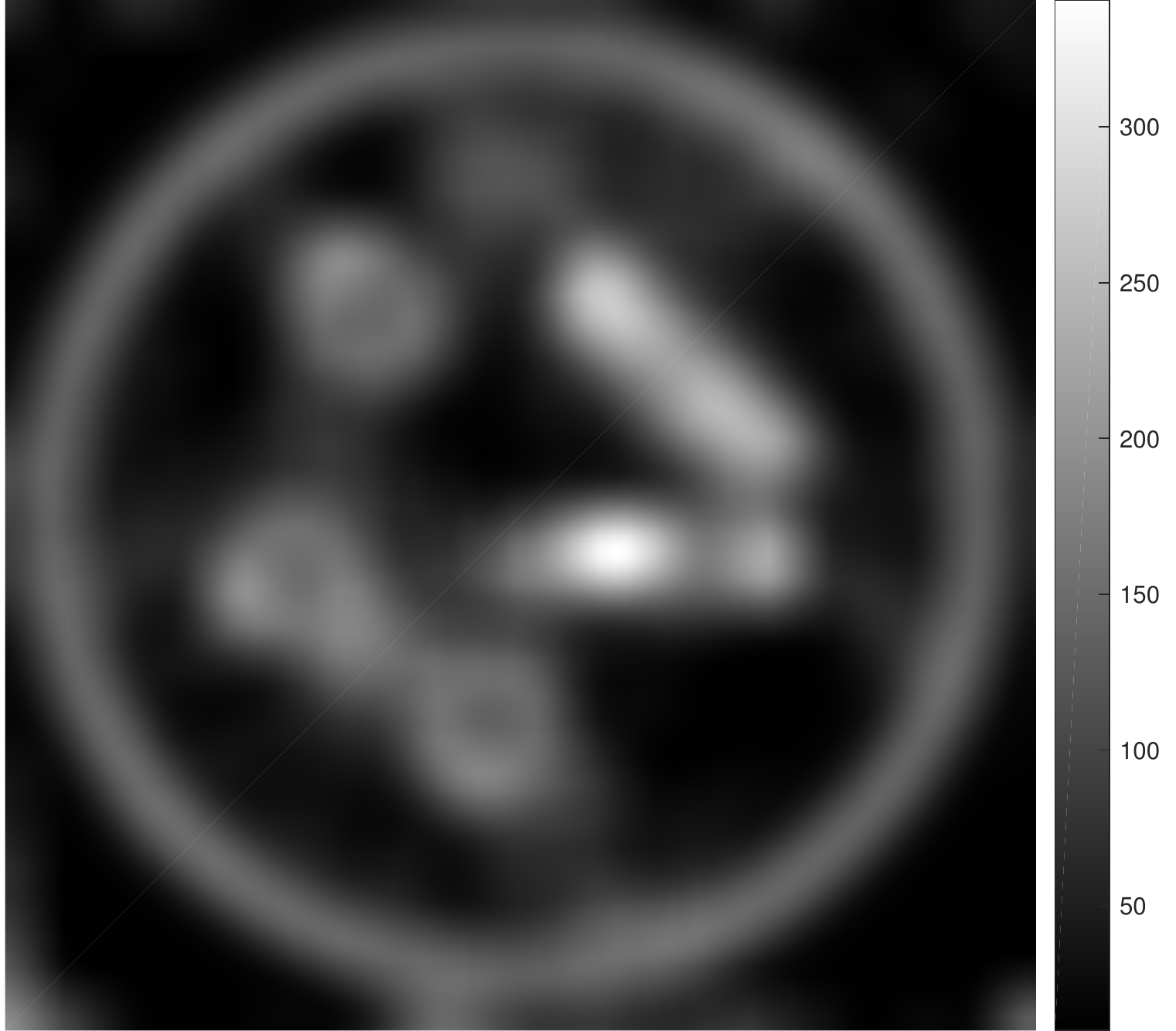}\\
\centering{(d)}
\end{minipage}
\end{center}
\caption{\label{result_gel2} Results of different methods for reconstructing the gel phantom with underdetermined rate 38\%. (a) FBP, (b) L2-TV with scalar $\lambda$, (c) Our method, (d) $\lambda$ in our method.}
\end{figure}

\textbf{Example 2.} Additionally to the simulated data, we also test our method on real CT measurements. In this experiment the gel phantom shown in Figure \ref{phantom} (b) is measured using a fan-beam geometry with 560 beams and 360 or 180 projection angles. The reconstructions are in a square domain of 512$\times$512 pixels, which result in an under-determined rate of 77\% and of 38\%, respectively. 

In Figure \ref{result_gel1} and \ref{result_gel2} we compare our method with the FBP algorithm and the L2-TV method. Due to insufficient measurements and noise, FBP cannot provide satisfactory results. Comparing the results obtained by our method with the ones from the L2-TV method, we see that our method reduces more artifacts while keeping similar quality on reconstructed object textures. Furthermore, from the final values of $\lambda$ obtained in our method we find that our method can correctly distinguish textured regions from homogeneous regions. Then, by setting different regularization parameter values, we vary the strength of the smoothness in the different regions.

\section{Conclusion}\label{sec:conclusion}

In this paper, we introduce a new tomographic reconstruction method with spatially varying regularization parameter. By introduce an auxiliary variable, the new approach extends the idea proposed in the SA-TV method to general inverse problems, where the data-fitting term and the regularization term fall in different domain. Numerically we have shown that the new approach can reduce the influence of the noise better as well as preserving more details comparing with the state-of-the-art. One limitation of the new approach is that an input of the ``noise'' level in the reconstruction is required due to the use of the SA-TV method. A better artifact estimation method will be further discovered in the future work.

\section{Acknowledgments}
We thank Samuli Siltanen and Alexander Meaney from the Industrial Mathematics Computed Tomography Laboratory at the University of Helsinki for creating the gel phantom and providing its CT data. YD acknowledges the support of the National Natural Science Foundation of China via Grant 11701388. CBS acknowledges support from the Leverhulme Trust project on Breaking the non-convexity barrier, EPSRC grant Nr. EP/M00483X/1, the EPSRC Centre Nr. EP/N014588/1, the RISE projects CHiPS and NoMADS, the Cantab Capital Institute for the Mathematics of Information and the Alan Turing Institute. Both authors would like to thank the Isaac Newton Institute for Mathematical Sciences, Cambridge, for support and hospitality during the programme ``Variational methods and effective algorithms for imaging and vision'' where work on this paper was undertaken. This programme was supported by EPSRC grant number EP/K032208/1.

\bibliographystyle{plain}

\end{document}